# EULER'S GRAPH WORLD
## - MORE CONJECTURES ON GRACEFULNESS BOUNDARIES-I

Suryaprakash Nagoji Rao*
Surusha5152@hotmail.com
March-April 2014

**ABSTRACT**. Euler graphs are characterized by the simple criterion that degree of each node is even. By restricting on the cycle types yet additional intrinsic properties of Euler graphs are unveiled. For example, regularity higher than degree two is impossible within the class $\varepsilon_i$ of Euler graphs with one type of cycles $C_n$, $n \equiv i \pmod 4$, $i=0,1,2,3$. Further, graphs in $\varepsilon_i$ are planar for $i=1,2,3$. In the light of new properties of Euler graphs more gracefulness boundaries are conjectured for subclasses of Euler graphs and where relevant extended for general class of graphs. In absence of general analytical results much of the published papers resort to proving an infinite class of graphs graceful or nongraceful. The purpose of this paper is not to give families of graphs graceful or not. Instead, based on the available information expected gracefulness boundaries are proposed which may guide where to look for graceful graphs or lead to characterizations. While the (Ringel,Kotzig,Rosa) Tree Conjecture continues to remain unsettled, the work reported here serves an update on the conjectures made in Rao Hebbare (1975,1981) and more conjectures subsequently made in Rao (1999,2000) based on embedding theorems and graceful algorithms for constructing graceful graphs from a graceful graph. It is hoped that these conjectures lead to analytical techniques for establishing gracefulness property. Further probe into Euler graphs with only two types of cycles and other combinations of cycles continues.

## 1. INTRODUCTION

The word 'graph' will mean a finite, undirected graph without loops and multiple edges. Unless otherwise stated a graph is connected. For terminology and notation not defined here we refer to Harary (1972), Mayeda (1972), Buckly (1987). A graph G is called a 'labeled graph' when each node u is assigned a label $\varphi(u)$ and each edge uv is assigned the label $\varphi(uv)=|\varphi(u)-\varphi(v)|$. In this case $\varphi$ is called a 'labeling' of G. Define $N(\varphi)=\{n \in \{0,1,..., q_0\}: \varphi(u)=n,$ for some $u \in V\}$, $E(\varphi)=\{e \in \{1,2,...,q_0\}: |\varphi(u)-\varphi(v)|=e,$ for some edge $uv \in E(G)\}$. Elements of $N(\varphi)$, $(E(\varphi))$ are called 'node (edge) labels' of G with respect to $\varphi$. A (p,q)-graph G is 'gracefully labeled' if there is a labeling $\varphi$ of G such that $N(\varphi) \subseteq \{0,1,...,q\}$ and $E(\varphi)=\{1,2,...,q\}$. Such a labeling is called a 'graceful labeling' of G. A 'graceful graph' can be gracefully labeled, otherwise it is a 'nongraceful graph' (see Rosa (1967), Golomb (1972), Sheppard (1976) and Guy (1977) for chronology 1969-'77 and Bermond (1978). If a graph G is nongraceful, then a labeling $\varphi$ which gives distinct edge labels such that the maximum of the node labels is minimum is called an 'Optimal labeling' of G and the graph is called 'optimally labeled graph'. Note that $0 \in N(\varphi)$. This minimum is denoted by $opt(G) \geq q$ with equality holding whenever G is graceful. A node u (an edge e=uv) of G is called i-attractive (i-repelling) if there (there does not) exist a graceful labeling $\alpha$ of G such that $\alpha(u)=i$ ($|\alpha(u)-\alpha(v)|=i$). A graph G is 'highly graceful' if every connected subgraph of it is graceful.

A *cycle graph* is a graph that consists of a single cycle. *Pendant free graph* G has node degrees two or more. Euler graphs are pendant free graphs. *Core graph* is obtained from G by deleting all pendant nodes recursively till no pendant nodes exist. Core graph of a graph is pendant free. Core graph of a tree is empty graph and core graph of a unicyclic graph is cycle graph. Core graph of pendant free graph is the graph itself.

By *planting a graph G onto a graph H* we mean identifying a node of G with a node of H. For a pendant free graph G of order p, a *graphforest* GF(G) is constructed from G by planting any number of trees at each node of G. GF(G) is trivial when it is of order p, ie., GF(G)=G and the trees planted are $K_1$s. We assume GF(G) is nontrivial and is of order larger than p. The smallest nontrivial graphforest has pendant node planted at a node. The choice of planting a tree at a node is random. That is, trees of any order and in any number may be planted at each node. The requirement of G a pendant free graph is for convenience. If G is not pendant free then start with core graph of G.

Note that for any graph G the simplest graph structure containing G of higher order is graphforest. A graphforest is *graphtree* when one nontrivial tree is planted on one node. Graphforest of Euler graph G is called *Eulerforest* EF(G). *Cycleforest* $CF(C_n)$ and Treeforest TF(T) are similarly defined. The class of cycleforests is precisely the class of unicyclic graphs. Note that graphforest of a forest is forest and graphforest of a tree is a tree. The terms like

---

Eulertree, cycletree may be similarly defined.

Varied applications of labeled graphs have been cited in the literature Gallion (2013), for example, communication Networks, X-Ray Crystallography, Radio Astronomy, Circuit Layout Design and Missile Guidance. Labeled directed graphs were studied and applied to Algebraic Systems, Generalized Complete Mappings, Network Addressing Problems and N-Queen Problems (Gallion (2013): Bloom (1977), Bloom and Golomb (1977,1978), Harary (1988)). An application in the domain of networks is by I.C. Arkut, R.C. Arkut and N. Ghani (2000) on positive effects of the new Multi Protocol Label Switching (MPLS) routing platform in IP networks.

The simplest class of connected graphs is trees. Even this class exhibits complexity with respect to graceful labeling. Attempts to prove them graceful lead to the following:

**Conjecture** A (Ringel (1964), Rosa (1967,1991), Kotzig (1973)). All trees are graceful.

If true implies that all trees are highly graceful. Further, truth of the conjecture implies that a nongraceful graph has a cycle. For a detailed exposition on this conjecture we refer to Ringel (1964), Kotzig (1973), Rosa (1967,1991), Golomb (1972), Bermond (1978), Bloom (1979). For a recent survey on graceful trees we refer to (Robeva (2011)) although there are papers posted on web claiming to have settled the tree conjecture.

Erdos proved the following:

**Theorem** A (Paul Erdos, See Golomb (1972)). Almost all graphs are nongraceful.

So it is clear that graceful graphs are rare in the class of graphs. Complete graphs are nongraceful for orders five or more. However, other extreme viz., complete bipartite graphs are graceful. Many infinite classes of bipartite graphs and trees are known graceful. Bipartite nongraceful graphs so far known are Rosa-Golumb type. Contrary to Theorem A, for bipartite graphs:

**Conjecture** B (Gangopadhyay and Rao Hebbare (1980)). Almost all bipartite graphs are graceful.

Cycle graphs were shown graceful by Rosa (1967). Cycle graphs and wheel graphs were considered for gracefulness in the thesis of Rao Hebbare (1975). Cycle graphs were given optimal labelings. Wheel graphs were conjectured graceful. Trees satisfy $p=q+1$ and so every graceful labeling of a tree uses all node labels $0,1,\ldots,q$. whereas unicyclic graphs satisfy $p=q$ and so every graceful labeling misses a unique node label. Graceful labelings for cycle graphs of orders $p=3,4,7,8,11,12$ were enumerated and classified using missing node label. In case of nongraceful graphs optimal labeling additionally has missing edge labels and for nongraceful cycle graphs missing edge label is unique. A catalogue of nongraceful graphs of orders 5 and 6 were given. R. Frucht (1979) showed that wheels are graceful.

A conjecture (Rao Hebbare (1975)) that graphs, not of Rosa-Golumb type, with at least two blocks, each block a complete graph with at least three nodes, are nongraceful is false. A counterexample is shown in Fig.1a. The above class may be restricted to graphs with two complete blocks. Windmill graphs $mK_n$ ($n \geq 3$) consist of m copies of $K_n$ with a node in common. Windmills with n=3 (n=4) are named Dutch (French) m-windmills. In general, asymmetric windmill consisting of m complete graphs with a common node may be defined. Windmill of cycles $mC_n$ consist of m copies of n-cycles with a common node. Graceful labelings for $mC_4$, for m=2,3,4,5 and optimal labelings for the nongraceful graphs $2K_3$, $3K_3$, $2K_4$, $K_4K_3$ (a graph with $K_4$ and $K_3$ with a node in common) were given. Snakes of cycles and/or complete graphs may be similarly defined. Triangular snake with three triangles is graceful. Another near complete graph has m copies of $K_n$ with $K_{n-1}$ in common. This for n=3, an infinite class of graphs was shown graceful.

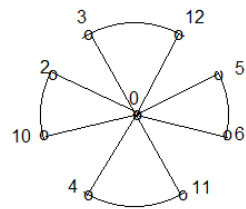

Fig.1a Graceful Euler (9,12)-graph.

**Theorem** (J.C. Bermond, A. Kotzig and J. Turgeon (1978)). The Dutch m-windmill is graceful if and only if $m \equiv 0$ or $1 \pmod 4$.

The following conjecture is still open. It was verified for $4 \leq m \leq 32$ (J. Huang and S. Skiena (1994)):

**Conjecture** (J.C. Bermond, A. Kotzig and J. Turgeon (1978)). The French m-windmill is graceful if m≥4.

Based on graceful algorithms and embedding theorems, conjectures were made in Rao Hebbare (1981) and are still open. Work accomplished was reported in a series of reports under the Research Project No. HCS/DST/409/76, Mehta Research Institute, Allahabad (1980-81). Author's several attempts in search of graceful graphs satisfying certain properties or finding a graceful labeling for a specific class of graphs or proving a class of graphs nongraceful lead to conjectures on the expected boundaries between gracefulness and non-gracefulness. With the advent of the mainframe computers and PCs while with ONGC, programs in FORTRAN or MATLAB were written for finding or enumerating graceful labelings thus proving a graph graceful or not. Improved versions of Rao Hebbare (1981) were reported or presented at the conferences Rao (1999,2000). The work of Rao (1999) accompanied with diagrams and figures but the paper was published without the diagrams and figures.

Order, size and structure of a graph dictate the complexity of finding a graceful labeling. Proving a graph graceful or not, draws a lot of computing resource and time. As the order and size of the graph increase proving the existence or nonexistence of a graceful labeling even by computers not only demands efficient algorithm but also a large amount of processor time. (See Gallion (2013): T.A. Redl (2003), B. Smith and K. Petrie (2003), B. Smith (2006), Jean-Francois Puget and Barbara M. Smith (2010), Michelle Edwards (2006), Eshgi (2004,2010)).

In this paper, expected gracefulness boundaries are proposed in the light of new and simple properties of Euler graphs with restrictions on the cycle structure. They are further extended for general class of graphs where relevant.

The simplest structure beyond a given graph is graphforest containing the graph as induced subgraph. In the event the graph G itself graceful then it is expected that a graphforest of G is graceful (See Fig.1b). Algorithms supporting this were given in Rao (1981,1999,2000). In the figures green and pink shades mean gracefulness and Non-gracefulness property respectively.

**Conjecture 1.** Graphforest GF(G) of a graceful pendant free graph G is graceful.

This conjecture if true implies that pendant free graph of G or core graph of a given graph is the minimum order to establish gracefulness of G and graphforests of G.

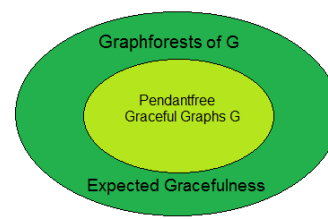

Fig.1b Expected gracefulness in Graphforests of pendantfree graphs.

## GRACEFUL ALGORITHMS

A *'graceful algorithm'* is one which when performed on a graceful or non-graceful graph as input gives rise to a graceful graph as output. In this section we review the results including four general existential graceful algorithms given in Rao Hebbare (1981) and Rao (1999). New graceful graphs were constructed from a given graceful or non-graceful graph exhibiting an optimal labeling. General and computationally efficient algorithms were given with illustrations. First we translate the algorithms in the form of theorems below:

**Theorem** (Rao (1999)). If G is graceful then a cater pillar planted on any 0-attractive node is graceful.

**Theorem** (Rao Hebbare (1981)). If G is a graceful (p,q)-graph G *then* there is a graceful (p+m,q')-graph with m>0 and q'>q.

$H(G,\varphi)$ denotes the set of all graceful graphs constructed from G, a graceful labeling $\varphi$, as in the proof of theorem above. Three important corollaries may be inferred. When G is tree and H is connected then H is also a tree and the algorithm is called *graceful tree algorithm*.

**Tree Corollary** (Rao Hebbare (1981)). If T is a tree of order p and $\varphi'$ is a graceful labeling of T and u,v∈N(T) with $\varphi'(u)=0$ and $\varphi'(v)=p-1$, then uv∈E(T) the graphs constructed are graceful.

**Unicyclic Corollary** (Rao Hebbare (1981)). If T is a tree of order p and $\varphi'$ is a graceful labeling of T then each graph H(u,v) for any uv nonedge in T added is a graceful unicyclic graph of order p.

It is well known that addition of any edge in a tree results in a unicyclic graph. Note that, every graceful labeling of a tree leads to a large class of graceful unicyclic graphs of the same order for every right choice of the available missing edges.

**(Bipartite) Corollary** (Rao Hebbare (1981)). If G is bipartite graceful (p,q)-graph G with the bipartition $V=A\cup B$ and $\varphi'$ is a graceful labeling of G such that $0\in A$ and $\varphi'(u)<\varphi'(v)$ holds for any $u\in A$ and $v\in B$ then the graphs in $H(G,\varphi)$ are graceful.

The graphs of $H(G,\varphi)$ in this case may be bipartite graceful graphs depending on the edges added. In addition, large classes of graphs of same order but other sizes result depending on the number of missing node and edge labels in a graceful labeling. These algorithms are the fundamental blocks underlying the gracefulness boundary conjectures.

**EMBEDDING GRACEFUL GRAPHS**

Here we review the results obtained in Rao Hebbare (1981).

**Theorem** (Bloom (1975) and Acharya (1981)). Every graph can be embedded into a graceful graph.

**Theorem** (Rao Hebbare (1981)). Every graph can be embedded into a graceful graph as an induced subgraph.

An optimal graceful embedding H of a graph G with respect to a given optimal labeling satisfies that G is an induced subgraph of a graceful graph H with minimum possible nodes.

**Theorem** (Rao Hebbare (1981)). A (p,q)-graph with an optimal labeling has an induced optimal graceful embedding.

An optimal labeling of a nongraceful graph additionally admits both missing node and edge labels. The graph may be embedded into a graceful graph using the missing labels. In order to make it induced missing edges are introduced not between existing nodes. New nodes add to the order of graph embedding. Minimum possible order of the embedding is the optimal order.

**Conjecture** (Rao Hebbare (1981)). Every embedding of optimal order is connected.

**EULER GRAPHS**

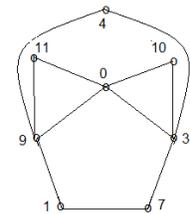

Fig.2 Example of Euler graph with 3,4,5,6,7-cycles.

Several properties of Euler graphs are known including characterizations (e.g., Zsolt (2010)). Cycle decomposition of Euler graph consists of cycles each edge occurring once. Euler graphs have no pendant nodes and their cycle decomposition may consist of cycles of any length. Fig.2 is a graph with cycles $C_n$, n=3,4,5,6,7. For a cycle decomposition $\xi_i$ denotes number of cycles of the type $C_n$, $n\equiv i(\mod 4)$, i=0,1,2,3. We shall focus on Euler graphs with restricted cycle structure or presence of certain types of cycles and so cycles in its cycle decompositions. Euler graphs exhibit yet more interesting properties so far unknown. As for example, regular Euler graphs of degree >2 with only cycles of the type $C_n$, $n\equiv 0(\mod 4)$ are nonexistent. These properties and others help in understanding boundaries of gracefulness leading to propose more conjectures.

Some characterizations and properties of Euler graphs (See Zsolt (2010)):

**Theorem** B. The number of edge disjoint paths between any two nodes of Euler graph is even.

**Theorem** C (Euler (1736), Hierholzer (1893)). A connected graph G is Euler graph if and only if all nodes of G are of even degree.

**Theorem** D (Veblen (1912-13)). A connected graph G is Euler if and only if its edge set can be decomposed into cycles.

**Theorem** E (Toida (1973), McKee (1984)). Fleischner (1989,1990).) A connected graph is Euler if and only if each of its edges lies on an odd number of cycles.

**Theorem** F (Bondy and Halberstam (1986)). A graph is Euler if and only if it has an odd number of cycle decompositions.

**Observation 1.** Every cycle decomposition of a bipartite Euler graph has only even cycles.

**Observation 2.** A non-bipartite Euler graph has a cycle decomposition consisting of an odd cycle. However, it may have cycle decompositions consisting of only even cycles. For example, Euler graph with a 4-cycle on each edge of a triangle is non-bipartite Euler graph and decomposes into three 4-cycles or into a 3-cycle and a 9-cycle.

Gracefulness exhibits affinity with Euler graphs as in:

**Theorem** G (Rosa (1967), Golomb (1972)). A necessary condition for Euler (p,q)-graph to be graceful is that $[(q+1)/2]$ is even.

Euler graph of this type is called here a *'Rosa-Golomb graph'*. This implies that Euler graphs with $q \equiv 1 \text{ or } 2 \pmod 4$ are nongraceful. Most of the known nongraceful graphs contain a subgraph isomorphic to Rosa-Golomb graph. The necessary condition for Euler graphs leads to understanding the intrinsic cycle structural properties of graceful Euler graphs. This as a guide first we shall study Euler graphs with only one type of cycles. We consider four subclasses of Euler graphs. Denote by $\varepsilon_i$ the class of Euler graphs with only cycles $C_n$, $n \equiv i \pmod 4$, $i = 0,1,2,3$. For example, $\varepsilon_0$ is the class of Euler graphs having only cycles of the type $C_n$, $n \equiv 0 \pmod 4$. $K_3$ is the smallest Euler graph and is graceful. Note that $K_n$ is nongraceful for $n > 4$ and Euler for $n$ odd. The smallest nongraceful Euler complete graphs and not Rosa-Golomb type graph are $K_9$ and $K_{11}$ with $q = 36 \equiv 0 \pmod 4$ and $q = 55 \equiv 3 \pmod 4$ respectively. They are regular Euler graphs of degree 8 and 10 respectively.

Consider the family of graphs $H(l,m,n) = mK_l \times \overline{K_n}$. A conjecture by Rao Hebbare (1975) on a simple infinite class of (n+4,4n+2)-graphs $G$ that $2K_2 \times \overline{K_n}$, $n > 1$, $n$ even are nongraceful was settled affirmatively. The graph $G$ is Euler whenever the degree of end nodes of $2K_2$ is even. This is the case when $n$ is odd. $q = 4(2t+1)+2 = 8t+6$ and so $q \equiv 2 \pmod 4$ and so $G$ is nongraceful. When $n$ is even $G$ is not Euler and $q \equiv 2 \pmod 4$, $G$ needs verification for graceful or nongraceful property. The proof technique uses a result by Golomb (1972) that in a graceful graph the sum of the edge weights on any cycle is even. Note that the graphs are not completely Rosa-Golomb type.

**Theorem** H (Bhat-Nayak and S.K. Gokhale (1986)). $H(2,2,n)$ is nongraceful, $n > 0$.

Some known nongraceful graphs enumerating graceful labelings on computer employing efficient algorithms based on Mathematical, CSP Model, Constrained Programming, Heuristic Programming are (see Gallion (2013): T.A. Redl (2003), B. Smith and K. Petrie (2003), B. Smith (2006), J.F. Puget and B.M. Smith (2010), M. Edwards (2006), Eshgi(2003)): $K_n \times P_2$ for $n = 6,7,8,9$; $K_n \times P_4$; $K_n \times C_3$, for $n = 3,5,6$; $K_n \times C_4$, for $n = 3,4$; Double Wheel: $DW_3$; $B(n,r,m)$ has $m$ cliques of size $n$ sharing a clique of size $r$: $B(n,3,m)$ for $n \geq 11$, $m \geq 2$; $B(6,2,2)$, $B(7,2,2)$, $B(8,2,2)$, $B(6,3,2)$, $B(6,4,2)$, $B(6,5,2)$, $B(6,3,3)$, $B(7,3,2)$.

## GRACEFULNESS BOUNDARIES THROUGH EULER GRAPHS

### Euler Graphs with Only One Type of Cycles

As an extension of the graceful tree conjecture we investigate graphs just next in complexity to the class of trees. It is well known that a graph which is not a tree has a cycle. Cycle graphs and unicyclic graphs are the simplest graphs with a unique cycle. Euler graphs admit only even degrees. The necessary condition for gracefulness for Euler graphs is a clue and guides for an extension of the class of trees towards gracefulness. With this goal first we shall prove some new simple properties of Euler graphs with restrictions on the cycle structure. A natural generalization of trees is the class of graphs with only one type of cycles. That is, cycles in any cycle decomposition of Euler graph are only one type. In the ensuing paragraphs Euler graphs with only one type of cycles are considered under the operation (mod 4). This simplifies the structure of Euler graphs and stand next to trees. It emerges that regularity of degree greater than two under this condition is impossible within the class of Euler graphs. That is, we prove in the ensuing paragraphs that regular Euler graphs with just one type of cycles are cycle graphs. Extremal graphs with the property are of interest. Some extremal Euler graphs with cycles $C_n$, n≡i(mod 4), i=0,1,2,3 are given.

**Observation 3**. Graphs with each block a cycle graph $C_n$, n≡i(mod 4) belongs to $\varepsilon_i$, i=0,1,2,3. In general, a graph G belongs to $\varepsilon_i$ iff each of its blocks belongs to $\varepsilon_i$, i=0,1,2,3.

In view of this, it is safe to assume that the graphs in the following discussion are blocks, that is, without cutnodes unless otherwise stated.

**Observation 4**. Suppose G is a graph from $\varepsilon_i$ and C a cycle in G. The graph H with the edges of C deleted from G belongs to $\varepsilon_i$ or each nontrivial component of H belongs to $\varepsilon_i$, i=0,1,2,3.

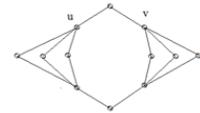
Fig.3a

Node pairs u,v each of degree ≥4 having at least four edge disjoint u-v paths are of interest. This may not hold for all node pairs follows from an example of Euler graph shown in Fig.3a with two nodes u,v of degree 4 but having no more than two edge disjoint paths. However, existence of such a pair is assured in any Euler graph as in:

**Observation 5**. For Euler graph G, if u is a node with d(u)≥4 then there is v, d(v)≥4 with at least four edge disjoint u-v paths. Clearly G is not a cycle graph. Node u is not unique otherwise it follows that u is a cutnode. So there are at least two nodes with degree ≥4 in G. Suppose u,v are such nodes. Since G is a block there is a common cycle C containing u,v. Two cases arise according as both u,v belong to same component or not in the graph H obtained by deleting the edges of C from G. In the first case, the component containing both u,v is Euler and u,v are of degree ≥2. The number of edge disjoint u-v paths in H is even and so there are at least four edge disjoint u-v paths in G. In the second case, u,v belong to different components in H. Consider the component containing u and rename the farthest node along C in the component by v and the result follows from u,v.

**Intersecting Cycles.** Two cycles may have nodes and paths in common. In general, the cycles may have multiple paths in common. Here, a path may be simply a node. By *two cycles intersect* we mean they have only one path in common with the nodes and edges of the path in common.

**Observation 6**. A graph from $\varepsilon_i$, i=0,1,2,3 which is not a cycle graph contains a series of cycles between u and v with a node in common with the next cycle. Consider a u-v path $P_1$ in black in a graph from $\varepsilon_i$, i=0,1,2,3, see Fig.3b. By Theorem B the number of edge disjoint u-v paths is even. A general second edge disjoint u-v path $P_2$: u=$u_0$,$u_1$;…;$u_{\alpha-1}$,$u_\alpha$=v; in red may have common nodes with $P_1$ in addition to u,v as shown. Suppose, $u_1$,$u_2$, …,$u_{\alpha-1}$ be the other common nodes. This results in a series of $\alpha>0$ cycles between nodes $u_i$ and $u_{i+1}$, i=0,1,…,$\alpha$-1. Any edge disjoint u-v path forms a series of cycles, two consecutive cycles having a node in common.

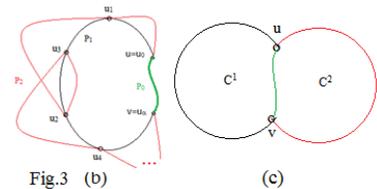
Fig.3 (b)    (c)

**Observation 7**. A graph from $\varepsilon_i$, i=0,1,2,3 which is not a cycle graph contains two edge disjoint cycles with only nodes u,v and a u-v path in common. Two u-v edge disjoint paths $P_0$ in green and $P_1$ in black of a graph in $\varepsilon_i$, i=0,1,2,3 form a cycle C. The node v may be so chosen so that there is no path between nodes of $P_0$ and $P_1$. $P_2$ be a third edge disjoint u-v path in red with no common nodes than u,v and such a choice is always possible. Then the union of the paths $P_0$, $P_1$ and $P_0$, $P_2$ form cycles with $P_0$ in common, see Fig.3c. Otherwise, $P_2$ has at least one common node with $P_0$, $P_1$ besides u,v. A general configuration of these paths with $P_2$ having no common node with $P_0$, is shown in Fig.3b. Let $u_1$ be the first common node from u between $P_1$ and $P_2$. Then there is a cycle containing u and $u_1$. If $u_1$, …,$u_\alpha$ are common nodes then there are $\alpha>0$ cycles containing u=$u_0$,$u_1$; …;$u_{\alpha-1}$,$u_\alpha$; each of these cycles $C_i$ satisfy ≡i(mod 4). The cycle $C^1$ containing u and $u_1$ and the cycle C have u-$u_1$ path in common as required. Any such edge disjoint u-v path forms a series of cycles any two consecutive cycles having a node in common. Further, each such cycle consists of two parts with a part of C and a part of $P_2$.

**Case-0. $\varepsilon_0$: Euler graphs with only cycles $C_n$, $n\equiv 0\pmod 4$**

**Construction**. If G is a (p,q)-graph in $\varepsilon_0$ then 1-subdivision graph $S_1(G)$ of G belongs to $\varepsilon_0$. $S_1(G)$ is Euler (p+q,2q)-graph and belongs to $\varepsilon_0$. Note that a cycle $C_n$ in G has n even and becomes $C_{2n}$ in $S_1(G)$. $K_{2,2n}$ for n>0 belongs to $\varepsilon_0$. In general, graphs of the type $\varepsilon_0$ may be constructed from G by adding any number of paths of same length between two nodes u,v at even distance making sure that every cycle in G containing such a path is of length $\equiv 0\pmod 4$. In particular, for any nonplanar bipartite graph G from $\varepsilon_0$, $S_1(G)$ is nonplanar and belongs to $\varepsilon_0$.

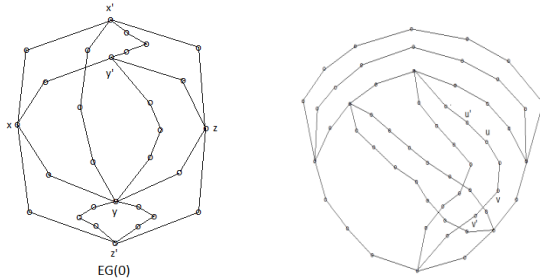

Fig.4a EG(0)

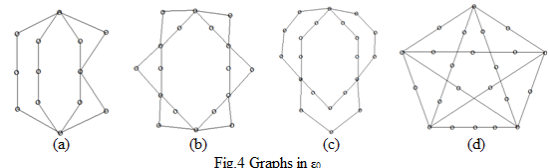

Fig.4 Graphs in $\varepsilon_0$

**Examples**: (i) $\varepsilon_0$ is a subclass of bipartite graphs containing only cycles of length multiples of 4. $K_{2,m}$ is the smallest Euler (m+2,2m)-graph of order m+2 from $\varepsilon_0$ for m even with node degrees 2 or m. All possible cycles are 4-cycles. This class of complete bipartite graphs is graceful. Three graphs with 4,8,12-cycle types from $\varepsilon_0$ are shown in Fig.4a,b,c. (ii) Subdivision graph may be generalized to include an arbitrary number of nodes newly introduced on each edge so that the resulting graph belongs to $\varepsilon_0$. An isolated example (23,28)-graph H homeomorphic to $K_5$ from $\varepsilon_0$ is shown in Fig.4d. $S_1(H)$ is a graph from $\varepsilon_0$. This gives rise to an infinite class of nonplanar graphs in $\varepsilon_0$. Its gracefulness is not known. Study the properties and characterize planar graphs of $\varepsilon_0$. (iii) Janakiraman and Sathiamoorthy (2013), have shown that Hanging Theta Graphs are graceful. These graphs have cycles of length $\equiv 0\pmod 4$ and so belong to $\varepsilon_0$.

**Theorem 2.** Every (p,q)-graph G in $\varepsilon_0$ satisfies that G is bipartite and $q\equiv 0\pmod 4$.

**Proof**. Every cycle decomposition of G consists of only cycles of the type $C_n$, $n\equiv 0\pmod 4$. Since each cycle is of length $n\equiv 0\pmod 4$ it follows that G is bipartite. Further, a cycle decomposition of G includes every edge exactly once and so $q\equiv 0\pmod 4$. □

**Theorem 3.** If two cycles of a graph in $\varepsilon_0$ intersect then the common path $P_t$ has even number of edges.

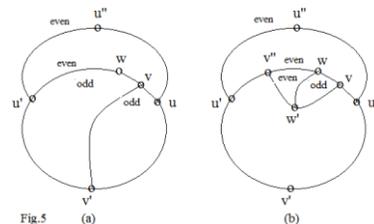

Fig.5  (a)   (b)

**Theorem 4.** A graph in $\varepsilon_0$ has at least one degree 2 node.

**Proof**. If G is a cycle graph then the result is trivially true. Otherwise, there is a node say u, of degree >2. Let u,v,w be a 2-path along a cycle. Then there are even numbers of u-v edge disjoint paths in G. Claim that v is a node of degree 2. If d(v)>2 then two cases arise, see Figs.5a,b. Firstly, let v' be such that v-v' is a path and v' is on a v-u'-u path along the cycle. But then the edge uv is the intersection of the cycles u,v-v'-u and u,v,w-u', a contradiction. The case v' is on u-u''-u' path is similar. Lastly, if v'' is a node on v,w-u' path then u'-v'' path is of even length as it is intersection of two cycles of $\equiv 0\pmod 4$. Similar to u and u', the conclusion for v and v'' holds, that is, for w there is w' on v-w'-v'' path, so that there is a w-w' path. But then wv is intersection of the two cycles w-w'-v,w and v,w-v''-u'-u''-u,v, a contradiction, completing the proof. □

In fact, we shall establish stronger result. For this we need a result for any adjacent nodes as follows:

**Theorem 5.** A graph in $\varepsilon_0$ satisfies that every pair of adjacent nodes has exactly two edge disjoint paths.

Proof. Suppose G is a graph from $\varepsilon_0$. If G is a cycle graph then the result is trivially true. Otherwise G has a node,

say u, of degree >2. For any two nodes u,v there are even number of edge disjoint u-v paths by Theorem B. In particular, this is true for any adjacent nodes. If d(v)=2 then the result holds. Otherwise, by Observation 5 there is a node, say v renaming if necessary, with d(v)>2 so that there are at least 4 edge disjoint u-v paths in G. The edge uv together with second edge disjoint u-v path results in a cycle, say $C^1$ of $\equiv 0 \pmod 4$. Consider a third edge disjoint u-v path. If it is node disjoint with $C^1$, except for u and v, then there are two edge disjoint cycles, as in Fig.3c, with the only edge uv in common. But this results in a combined cycle of $\equiv 2 \pmod 4$ with the edge uv deleted, a contradiction. Otherwise, every third edge disjoint u-v path has some nodes in common with $C^1$ besides u,v. Any such edge disjoint u-v path forms consecutive cycles with a node in common as in Observation 7. Further, each such cycle consists of a part of $C^1$ and a part of third edge disjoint u-v path. Both parts of such a cycle are of even length as all cycles are of $\equiv 0 \pmod 4$ in G. Therefore, it follows that length of every edge disjoint u-v path is of even length. However, by selection $C^1$ is $\equiv 0 \pmod 4$ and both u-v paths along $C^1$ are odd, a contradiction. We conclude that there are exactly two edge disjoint u-v paths for any adjacent nodes u,v in G. □

The result may be extended as follows:

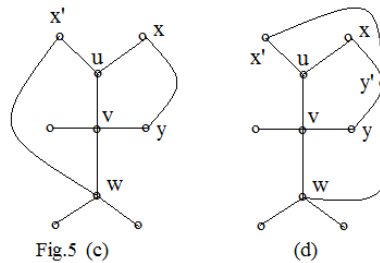

**Theorem 6.** A graph in $\varepsilon_0$ satisfies that every pair of nodes at odd distance, without edge disjoint u-v-w path with v in common, has exactly two edge disjoint paths.

Proof follows on similar lines of Theorem 4. We shall now deduce a stronger result than Theorem 4:

Fig.5 (c)   (d)

**Corollary 6.1.** At least one of three consecutive nodes on any path of a graph in $\varepsilon_0$ is of degree two.

**Proof.** Suppose u,v,w nodes form a u-w path with each node of degree >2. Let C be a cycle containing u,v. By Theorem 4 the edge uv has one more edge disjoint u-v path namely, u,x-y,v. If w is on this cycle then choose w adjacent to v not on this cycle containing u,v as shown in Fig.5c. Such a choice is assured as d(v)>2. Now consider x' adjacent to u. for x',u,v,w there is another edge disjoint x'-w path. If this is node disjoint with x',u,v,w except for x',w then u,v has a third edge disjoint u-v path, a contradiction. Otherwise there is a node say y' in common with u,x-y,v as shown in Fig.5d. But then this gives rise to a third edge disjoint u-v path, u-x'-y'-w,v, a contradiction implying that one of u,v,w is of degree 2. This completes the proof. □

An important consequence of Theorem 4 or Corollary 6.1:

**Corollary 6.2.** A graph in $\varepsilon_0$ is regular iff it is the cycle graph $C_n$, $n \equiv 0 \pmod 4$.
                    Or
         Regular Euler graphs of degree >2 in $\varepsilon_0$ are nonexistent.

**Corollary 6.3.** Regular bipartite Euler graph of degree >2 has at least one cycle of both types $C_n$, $n \equiv 0 \& 2 \pmod 4$.

**Theorem 7.** Size of a graph G $\in \varepsilon_0$ satisfies $q \equiv 0 \pmod 4$.

The size q of G satisfies: $q = 4t_0 * \xi_0$, $t_0 > 0$. In other words, for any $\xi_0 > 0$, $q \equiv 0 \pmod 4$. Therefore, graphs in $\varepsilon_0$ are not Rosa-Golomb type and so are candidates for (highly) gracefulness. Generalized graceful tree conjecture within the class of Euler graphs is:

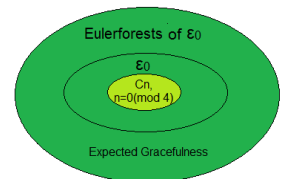

**Conjecture C.** (Rao Hebbare (1981)). Every graph in $\varepsilon_0$ is highly graceful.

This for Eulerforest may be stated as (See Fig.5e):

Fig. 5e Expected gracefulness of Eulerforests of $\varepsilon_0$.

**Conjecture 2.** A graph G in $\varepsilon_0$ is graceful iff Eulerforest EF(G) is graceful.

In other words, a graph with cycles $C_n$, $n \equiv 0 \pmod 4$ is graceful. Note that this class of Eulerforests includes unicyclic graphs with cycle graphs $C_n$, $n \equiv 0 \pmod 4$ as core graph.

### Case-1. $\varepsilon_1$: Euler Graphs with only cycles $C_n$, $n\equiv 1\pmod 4$

The smallest Euler graph in $\varepsilon_1$ is $C_5$ and is unique. Higher order graphs will have each block a cycle $C_n$, $n\equiv 1\pmod 4$. Two larger order Euler graphs in $\varepsilon_1$ are (9,10)-graph with two $C_5$s with a node in common and the cycle graph $C_9$.

If G is a cycle graph $C_n$, $n\equiv 1\pmod 4$ then G is Rosa-Golomb graph and so nongraceful. The structure of graphs in the class $\varepsilon_1$ is simple (see Fig.6) as described in:

**Theorem 8.** Each block of a graph in $\varepsilon_1$ is a cycle $C_n$, $n\equiv 1\pmod 4$.

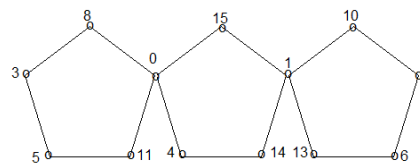

Fig. 6 Graceful graph from $\varepsilon_1$ with $q\equiv 1\pmod 4$.

**Proof.** By Observation 3 each block of a graph in $\varepsilon_1$ is also a graph in $\varepsilon_1$. Assume the graph is a block and not a cycle graph $C_n$, $n\equiv 1\pmod 4$ then there is a cycle of the same type with a path between two nonadjacent nodes of the cycle. This divides the cycle into two cycles with the path in common. But then one of the cycles is odd and the other is even, a contradiction. □

**Corollary 8.1.** Graphs in $\varepsilon_1$ are planar.

**Corollary 8.2.** A necessary condition for a graph in $\varepsilon_1$ graceful is that it has $4(t+1)$ or $4t+3$, $t\geq 0$ blocks.

**Corollary 8.3.** A graph in $\varepsilon_1$ is regular iff it is the cycle graph $C_n$, $n\equiv 1\pmod 4$.
    Or
Regular graphs of degree >2 in $\varepsilon_1$ are nonexistent.

**Theorem 9.** Size of a graph $G \in \varepsilon_1$ satisfies $q\equiv \xi_1\pmod 4$.

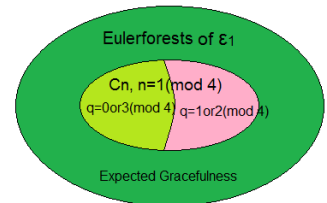

Fig. 7 Expected gracefulness of Eulerforests of $\varepsilon_1$.

The size of G satisfies: $q=(4t_1+1)*\xi_1$, $t_1>0$. In other words, for any $\xi_1>0$, $q\equiv \xi_1\pmod 4$. So if $\xi_1\equiv 0$ or $3\pmod 4$ then the graphs are candidates for gracefulness. If $\xi_1\equiv 1$ or $2\pmod 4$ then the graphs are nongraceful by Rosa-Golumb criterion.

**Conjecture 3.** A graph from $\varepsilon_1$ with $\xi_1\equiv 0$ or $3\pmod 4$ blocks is graceful.

Such (p,q)-graphs in $\varepsilon_1$ satisfy $q\equiv 0$ or $3\pmod 4$ and so candidates for gracefulness. For Eulerforests of $\varepsilon_1$ (See Fig.7):

**Conjecture 4.** For any graph G in $\varepsilon_1$, graceful or not, every Eulerforest EF(G) is graceful.

Note that this class of Eulerforests includes unicyclic graphs with cycle graphs $C_n$, $n\equiv 1\pmod 4$ as core graphs.

### Case-2. $\varepsilon_2$: Euler Graphs with only cycles $C_n$, $n\equiv 2\pmod 4$

This is a subclass of bipartite graphs. The smallest Euler graph in $\varepsilon_2$ is $C_6$. Next graph, not Euler, with cycles of $\equiv 2\pmod 4$ only is of order 8 with a path $P_3$ added between any two nodes of $C_6$ at distance 3. Another path $P_3$ added between the same pair of nodes results in Euler (10,12)-graph in $\varepsilon_2$.

**Construction.** If G is a (p,q)-graph in $\varepsilon_2$ then 2-subdivision graph $S_2(G)$ of G belongs to $\varepsilon_2$. In fact, if G is any Euler (p,q)-graph then $S_2(G)$ is Euler (p+q,3q)-graph and belongs to $\varepsilon_2$. The size of G satisfies $q\equiv 2\pmod 4$ and so $3q\equiv 2\pmod 4$. Further, a cycle $C_n$ in G becomes $C_{3n}$ in $S_2(G)$. In general, graphs of the type $\varepsilon_2$ may be constructed from G by adding any number of paths of same length between any two nodes u,v at odd distance making sure that every cycle in G containing such a path is $\equiv 2\pmod 4$.

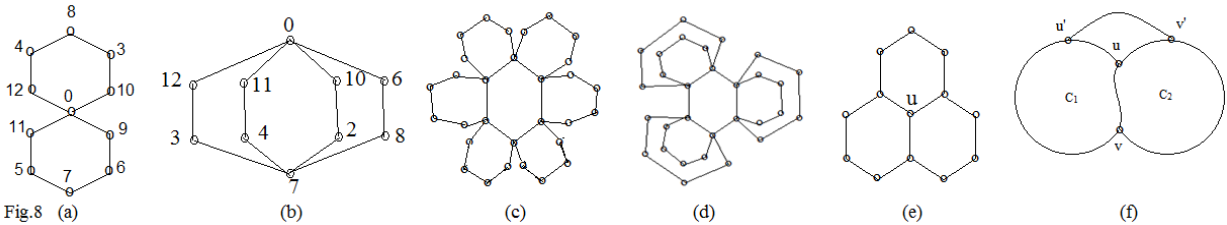

Fig.8 (a) (b) (c) (d) (e) (f)

**Examples**: Graceful graphs from $\varepsilon_2$ with p=10,11, q=12 are shown in Fig.8a,b. Two more graphs from $\varepsilon_2$ are: Fig.8c with a pair of nodes each of degree 4 and four edge disjoint paths and Fig.8d with a pair of nodes each of degree 4 but with exactly two edge disjoint paths.

**Observation 8.** The graph in Fig.8e is forbidden as subgraph for any graph in $\varepsilon_2$. Note that, there are three 6-cycles with each pair intersecting in an edge and two consecutive 6-cycles result in a combined 10-cycle if the common edge is deleted, but a 12-cycle results if the central node u is deleted.

**Observation 9.** For a node u' on $C^1$ and v' on $C^2$ there exists no u'-v' path for any graph in $\varepsilon_2$, see Fig.8f. This follows as any such path forms a cycle with one part even in common with $C^1$ or $C^2$, a contradiction.

**Theorem 10**. Every (p,q)-graph in $\varepsilon_2$ is bipartite with q≡0 or 2(mod 4).

**Proof**. Cycle decomposition of G in $\varepsilon_2$ has only even cycle types $C_n$, n≡2(mod 4). That the graph is bipartite follows. If $\xi_2$ is even (odd) then q≡0(mod 4) (q≡2(mod 4)).

When q≡2(mod 4), G is Rosa-Golomb and nongraceful and when q≡0(mod 4) the graphs are candidates for gracefulness and not highly graceful.

**Theorem 11.** If two cycles of a graph in $\varepsilon_2$ intersect then the common path $P_t$ has odd number of edges.

**Theorem 12**. Every pair of nodes u,w at distance 2 without edge disjoint u-v-w path with v in common in a graph G from $\varepsilon_2$ has exactly two edge disjoint paths.

**Proof**. Suppose, u,v,w is a 2-path in G with degree of u,w ≥4. Since G is Euler, by Theorem B even number of edge disjoint u-w paths is even. The first u-w path: u,v,w together with second u-w edge disjoint path results in a cycle, say $C^1$ of type ≡2(mod 4). Consider a third edge disjoint u-w path. If it is node disjoint except for u,v,w then there are two edge disjoint cycles each of type ≡2(mod 4), as in Fig.3c, with the only path u,v,w in common. But this results in a combined cycle of ≡0(mod 4) with the edges uv,vw and the node v deleted, a contradiction.

Otherwise, the third edge disjoint u-w path has at least one node in common with $C^1$. Let u' be the first node on $C^1$ with the edge disjoint u-u' path as shown in Fig.9a together with the u,v,w-u' path along $C^1$ forms a cycle and so is of ≡0(mod 4). Each of the two parts of any such cycle is odd. The node w has a third path say to w'. Two cases arise according as w' is on u'-u path or w-u' path as in Fig.9b in the first case u,v,w path is common between the cycles u,v,w-w'-u and u,v,w-u'u''-u, a contradiction. In the second case let x,y nodes be such that wx, xy are edges. The only case remains to be verified is when there is a y-y' path for y' on w-w' path not along $C^1$. Here the w,x,y path is common between the cycles w,x,y-y'-w and u,v,w,x,y-w''-u'u''-u, a contradiction. This completes the proof. □

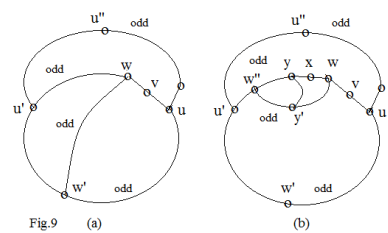

Fig.9 (a) (b)

This may be extended to:

**Theorem 13**. Every pair of nodes u,w at even distance without edge disjoint u-w path having at least one node of u-w path in common in a graph G from $\varepsilon_2$ has exactly two edge disjoint paths.

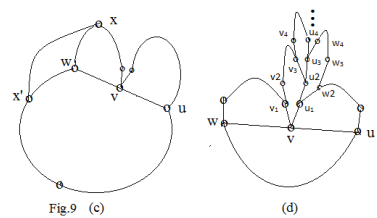

Fig.9 (c) (d)

**Theorem 14.** A graph in $\varepsilon_2$ has at least one degree 2 node.

**Proof**. Suppose G is a graph from $\varepsilon_2$ with no degree 2 node. Let u,v,w be a 2-path along a cycle, see Fig.9a. Then there are even, in this case at least four, numbers of edge disjoint u-w paths in G. Such a choice is possible follows from Observation 5. Two cases arise according as there is an edge disjoint u-w path with node v in common or not.

First suppose u,v,w and u-u''-u'-w be two edge disjoint u-w paths. Claim that w is a node of degree 2. If d(w)>2 then two cases arise. Firstly, let w' be such that w-w' is a path and w' is on a w-u'-u path along the cycle, see Fig.9a. But then the path u,v,w is the intersection of the cycles u,v,w-w'-u and u,v,w-u'-u''-u, a contradiction. The case w' on u-u''-u' path is similar. Lastly, if w'' is a node on w,x,y-u' path then w,x,y path is intersection of two cycles of $\equiv 0 \pmod 4$ as shown in Fig.9b. Similar to u and u', the conclusion for y and y' holds, that is, for y there is y' on w-y'-w' path, so that there is a y-y' path. But then w,x,y is intersection of two cycles w,x,y-y',w and u,v,w,x,y-w''-u'-u''-u, a contradiction.

Suppose u,$w_2$,$u_1$,v,$v_1$,w be an edge disjoint u-v path with only v in common, see Fig.9c. Firstly, no node, say x, of this path has a path to a node, say y, of u-x-w path as shown. If there is such a pair then the arcs x'w and wx are odd as they are intersection of cycles. Thus the arc x-w-x' is of even length. But the arc x'wx is intersection of the cycles x-x'-w-x and x-w-x'-uv-x and is odd, a contradiction. Lastly, neither of $u_1$ and $v_1$ is of degree 2 node, see Fig.9d. Suppose $u_2$ and $v_2$ are adjacent respectively. It follows that they are new nodes as $u_1$ and $v_1$ cannot be adjacent to any node on the cycle $C_1$ and u-$w_2$,$u_2$,v,u and v,$v_1$-w,v. By the same reason there is an edge disjoint $u_2$-$v_2$ path as shown. Next the choice of $u_2$, $w_2$ nodes lead to $w_3$,$u_3$ nodes as shown. The procedure may be continued to infinity with v,$u_1$,$u_2$,$u_3$, … path being infinite, contrary to finite graph. □

As an important consequence, we have

**Corollary** 14.1. A graph in $\varepsilon_2$ is regular iff it is the cycle graph $C_n$, n$\equiv$2(mod 4).
                         Or
         Regular graphs of degree >2 in $\varepsilon_2$ are nonexistent.

**Corollary** 14.2. Regular bipartite Euler graph of degree >2 with only cycles of the type $C_n$, n$\equiv$2(mod 4) also has at least one cycle of the type $C_n$, n$\equiv$0(mod 4).

Sum of two numbers of the type $\equiv$2(mod 4) may lead to one of the types $\equiv$0(mod 4) or $\equiv$2(mod 4). As a result the structural complexity of the graphs in $\varepsilon_2$ is simple in comparison with graphs in $\varepsilon_0$. In fact, we shall prove:

**Theorem 15.** Graphs in $\varepsilon_2$ are planar.

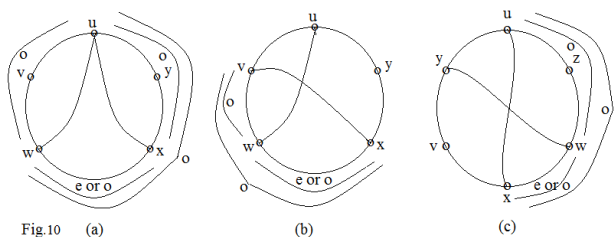

**Proof**. The graphs in Fig.10 are forbidden as subgraphs in a graph G from $\varepsilon_2$. Figs.10a,b are subgraphs of a graph homeomorphic to $K_5$ and Fig.10c is a subgraph of a graph homeomorphic to $K_{3,3}$. That the subgraphs are forbidden follows as arc wx can be neither even nor odd. The arcs may be viewed as intersection of cycles of the type $\equiv$2(mod 4) and they are of odd length by Theorem 11. □

**Theorem 16**. Size of a graph G $\in \varepsilon_2$ satisfies q$\equiv$ 2$\xi_2$(mod 4).

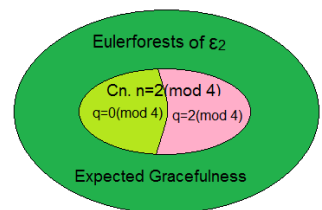

The size of G satisfies: q= (4$t_2$+2)*$\xi_2$, $t_2$>0. In other words, for any $\xi_2$>0, q$\equiv$ 2$\xi_2$(mod 4). So if 2$\xi_2$$\equiv$0(mod 4) or $\xi_2$ is even then the graphs are candidates for gracefulness. If 2$\xi_2$$\equiv$2(mod 4) or $\xi_2$ is odd then the graphs are nongraceful.

**Conjecture 5**. A (p,q)-graph in $\varepsilon_2$ with $\xi_2$ even is graceful.

For Eulerforests in $\varepsilon_2$ (See Fig.10d):

Fig.10d Expected gracefulness for Eulerforests of $\varepsilon_2$.

**Conjecture 6**. For any graph G in $\varepsilon_2$ Eulerforest EF(G) is graceful.

Note that this class of Eulerforests includes unicyclic graphs with cycle graph $C_n$, n≡2(mod 4) as core graph which is not graceful.

**Case-3. $\varepsilon_3$: Euler Graphs with only cycles $C_n$, n≡3(mod 4)**

$K_3$ is in $\varepsilon_3$. Graphs with q≡3(mod 4) are graphs with each block a cycle graph $C_n$, n≡3(mod 4). Two $K_3$s with a node in common is the next Euler (5,6)-graph in $\varepsilon_3$. Euler (7,9)-graph has three $K_3$s in chain. Euler (9,10)-graph consists of a $K_3$ and cycle graph $C_7$ with a node in common.

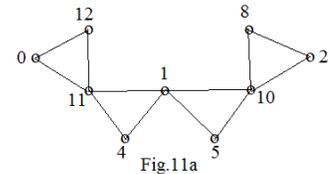

Consider the class of graphs obtained from a path $P_n$ or a cycle $C_n$ replacing each edge by a triangle. In other words, a graph is consisting of triangles either as a chain or a cycle. Such a graph is Euler (2n+1,3n)- or (2n,3n)-graph. The graph corresponding to a path is Euler and belongs to $\varepsilon_3$. For n=4 the graph is graceful shown in Fig.11a.
Examples: Graphs with each block a $C_n$, n≡3(mod 4) and number of blocks is 4t+1, t≥0.

The structure of graphs in this class is simple as described in:

**Theorem 17.** Each block of a graph in $\varepsilon_3$ is a cycle $C_n$, n≡3(mod 4).

**Proof**. By Observation 3 each block of a graph in $\varepsilon_3$ is also a graph in $\varepsilon_3$. Assume the graph is a block and not a cycle graph $C_n$, n≡3(mod 4) then there is a cycle of the same type with a path between two nonadjacent nodes of the cycle. This divides the cycle into two cycles with the path in common. One of the cycles is odd and the other is even cycle, a contradiction. □

**Corollary 17.1.** Graphs in $\varepsilon_3$ are planar.

Note that the cycle graphs $C_n$, n≡3(mod 4) are known to be graceful. Euler Graph with only cycles $C_n$, n≡3(mod 4) and no cutnodes is the cycle graph $C_n$, n≡3(mod 4). In general, a graph in $\varepsilon_3$ with a cutnode satisfies that each of its blocks is a cycle $C_n$, n≡3(mod 4). If the number of blocks in G is ≡2or3(mod 4) then q≡2or1(mod 4) respectively, so G is nongraceful else is a candidate for gracefulness.

**Corollary 17.2.** A graph from $\varepsilon_3$ is regular iff it is the cycle graph $C_n$, n≡3(mod 4).
                             Or
Regular graphs of degree >2 in $\varepsilon_3$ are nonexistent.

**Corollary 17.3.** A necessary condition for a graph in $\varepsilon_3$ graceful is that it has 4(t+1) or 4t+1, t≥0 blocks.

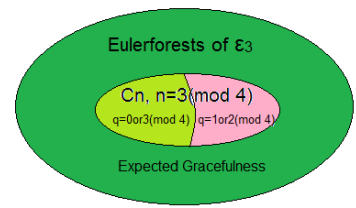

**Theorem 18**. Size of a graph G ∈ $\varepsilon_3$ satisfies q≡ 3$\xi_3$(mod 4).

Fig.11b Expected gracefulness of Eulerforests of $\varepsilon_3$.

The size of G satisfies: q=(4$t_3$+3)*$\xi_3$, $t_3$ >0. In other words, for any $\xi_3$>0, q≡3$\xi_3$(mod 4). So, if 3$\xi_3$≡0or3(mod 4) then the graphs are candidates for gracefulness. In this case, $\xi_3$≡0or1(mod 4). If 3$\xi_3$≡1or2(mod 4) then the graphs are nongraceful. In this case, $\xi_3$≡2or3(mod 4).

**Conjecture 7.** A graph from $\varepsilon_3$ with $\xi_3$≡0or1(mod 4) blocks is graceful.

Such a (p,q)-graph in $\varepsilon_3$ satisfies that q≡0or3(mod 4) and so are candidates for gracefulness. For Eulerforests of $\varepsilon_3$ (See Fig.11b).

**Conjecture 8.** For any graph G in $\varepsilon_3$ Eulerforest EF(G) is graceful.

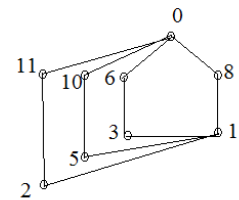

Fig.12a A graceful graph with cycles of length 5 or 6.

Note that this class of Eulerforests includes unicyclic graphs with cycle graphs $C_n$, n≡3(mod 4) as core graph.

A summary of regularity in the above four cases follows:

**Theorem 19.** Euler graph with only one type of cycles is regular iff it is regular of degree two and is the cycle graph.

**Euler Graphs with Mixed Cycles**

The graphs from $\varepsilon_0$ or $\varepsilon_3$ are friendlier to yield graceful labelings. Whereas graphs from $\varepsilon_1$ or $\varepsilon_2$ are either nongraceful or the graceful labelings become rare and difficult to find one. For example, graceful labeling of $K_5 \times P_2$ is unique not including complementary labeling. When mixed cycles are present even a combination of cycles ≡0or3(mod 4) may lead to a nongraceful graph. For example, the graph obtained from two $K_4$s with a node in common is known to be nongraceful. This graph has only cycles of length 3 and 4, but in combination still it exhibits nongracefulness. A general class of graphs of interest to probe gracefulness is: Two complete graphs $K_m$, $K_n$ with a node in common which also may be defined as Cartesian product of union of $K_{m-1}$ and $K_{n-1}$ with $K_1$, that is, $(K_{m-1} \cup K_{n-1}) \times K_1$, an $(m+n-1, {}^mC_2 + {}^nC_2)$-graph. An example of a graceful (9,11)-graph with cycles of length 5 or 6 as in Fig.12a.

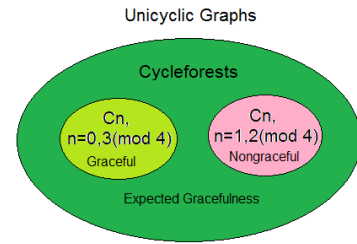
Fig.13a Expected gracefulness for Cycleforests.

It is expected that planting a tree structure onto a node of a graceful Euler graph shall not disturb gracefulness. Based on this intuition we generalize necessary part of Conjecture 2:

**Conjecture 9.** Eulerforest EF(G) of a graceful Euler graph G is graceful.

Cycleforests are unicyclic graphs. We state two conjectures for cycleforests which include unicyclic graphs and graphs with Rosa-Golomb graphs forbidden as subgraphs (See Figs.13a,b):

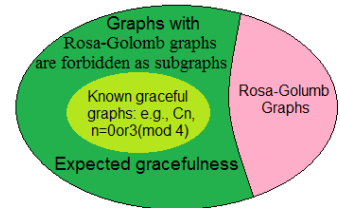
Fig.13b Expected gracefulness for graphs with Rosa-Golomb graphs forbidden as subgraphs.

**Conjecture D.** (Rao Hebbare (1981), Truszczynski (1984)). The only nongraceful unicyclic graphs are cycle graphs $C_n$, n≡1or2(mod 4).

**Conjecture E.** (Rao Hebbare (1981)) Euler graph G for which Rosa-Golomb graphs are forbidden is highly graceful.

If proved true implies that a bipartite Euler graph G without cycles $C_n$, n≡2(mod 4) is highly graceful. This essentially means that no other type of nongraceful bipartite Euler graphs is possible. In fact, the author's strong conviction is that highly graceful graphs are characterized by the simple criterion that Rosa-Golomb graphs are forbidden. The following extends Conjecture E:

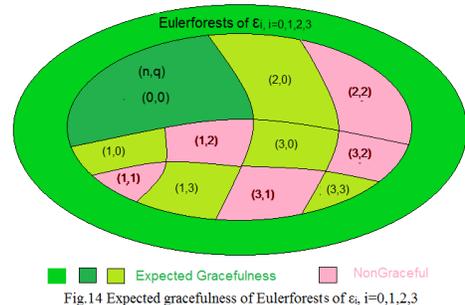
Fig.14 Expected gracefulness of Eulerforests of $\varepsilon_i$, i=0,1,2,3

**Conjecture F.** (Rao (1999)) A graph G with Rosa-Golomb graphs forbidden as subgraphs is highly graceful.

Expected gracefulness in the class of Euler graphs with only one type of cycles is summarized in Fig.14 The pair (i,j), i,j=0,1,2,3 denotes the class of Euler (p,q)-Graphs, with only type of cycles $C_n$, n≡i(mod 4) and q≡j(mod 4). Further probe into Euler graphs with combination of cycles continues.

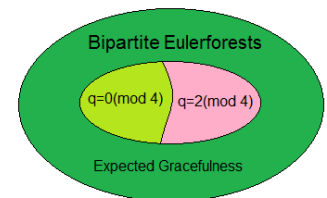
Fig.15 Expected gracefulness in bipartite Eulerforests

**BIPARTITE GRAPHS**

First we prove a result which helps intuitive process.

**Theorem 20**. Every bipartite graph G is embeddable as a subgraph in a bipartite Euler graph.

**Proof**. If G is Euler then we are done. Otherwise G has even number of odd nodes. If all odd nodes are in one partition then add a node in the other partition adjacent to all of the odd nodes. The resulting graph is Euler with G as a subgraph. Consider the case, odd nodes are in both the partitions in which case the number of odd nodes in each partition is either odd or even. If there are even number of odd nodes in each partition then the graph G' by adding a node in each partition and adjacent to the odd nodes of other partition is Euler with G as a subgraph. The only other case is that odd no of odd nodes in either partition. Add nodes in each partition adjacent to the odd nodes in other partition. Now the only odd nodes are the new nodes. Join the two odd nodes in each partition resulting in Euler with G as a subgraph. This completes the proof. □

Maximum node label in a labeling of a graceful or nongraceful graph satisfies that opt(G)≥q. This may get restricted in the class of bipartite Euler graphs as follows (See Fig.15):

**Conjecture** 10.   Bipartite Euler graph G satisfies that opt(G)=q or q+1.
                              Or
         Bipartite Rosa-Golomb (p,q)-graphs G are the only nongraceful bipartite Euler graphs
         and opt(G)=q+1 holds.

## Acknowledgements

The author, an Oil & Gas Professional with works in Graph Theory, expresses deep felt gratitude to
The *Sonangol Pesquisa e Produção*, Luanda, Angola
for excellent facilities, support and encouragement (2011-'13) with special thanks to
Mr. *Joao Noguiera*, Managing Director, Development Directorate.

Deep felt gratitude to *Alexander Rosa*, Emeritus Professor,
Dept. of Mathematics & Statistics, McMaster University, West Hamilton, Ontario, Canada and
*Solomon W Golumb*, an American Mathematician, Engineer and Professor of Electrical Engineering
at the University of Southern California (A recipient of the USC Presidential Medallion,
the IEEE Shannon Award of the Information Theory Society and three honorary doctorate degrees,
National Medal of Science presented by President Barack Obama)
for their initial inspiring works on Graceful Labelings and to
Emeritus Professor *G.A. Patwardhan*, Combinatorics, Department of Mathematics, IIT Bombay.